%% file: Intersection.tex
\documentclass[11pt,a4paper]{amsart}

\usepackage{amscd,amssymb,amsopn,amsmath,amsthm,graphics,amsfonts,mathrsfs,accents,enumerate,verbatim,calc}
\usepackage{amsmath}
\usepackage{amsthm}
\usepackage{amsfonts}
\usepackage{amssymb}
\usepackage{enumitem}
\usepackage{graphicx}
\usepackage[colorlinks=true,linkcolor=red,citecolor=blue]{hyperref}
\usepackage{latexsym}
\usepackage{mathrsfs}
\usepackage{makecell}  
\usepackage{newtxmath}
\usepackage{newclude}
\usepackage{tikz}		
\usepackage{tikz-cd}	
\usepackage[all]{xy}
\usepackage{nicematrix}
\usepackage{arydshln} 	
\usepackage{booktabs}
\usepackage{stmaryrd}
\usepackage[all]{xy}
\usepackage{float}

\date{}
\hypersetup
{colorlinks=true,linkcolor=black}
\allowdisplaybreaks[4] \footskip=15pt
\renewcommand{\uppercasenonmath}[1]{}

\numberwithin{equation}{section} \theoremstyle{plain}
\newtheorem{lemm}{Lemma}[section]
\newtheorem{coro}[lemm]{Corollary}

\newtheorem{theo}[lemm]{Theorem}
\newtheorem{remark}[lemm]{Remark}

\newtheorem{exam}[lemm]{Example}

\definecolor{morpink}{RGB}{239,229,248}
\definecolor{morblue}{RGB}{146,172,209}
\definecolor{defcolor}{RGB}{50,205,50}
\definecolor{tcolor}{RGB}{210,105,30} 
\definecolor{lcolor}{RGB}{210,105,30} 
\definecolor{pcolor}{RGB}{253,249,238} 

\def\no{\noindent}

\def\ra{\rightarrow}
\def\dan{\rightarrowtail}
\def\man{\twoheadrightarrow}

\def\proof{\no\rmm{Proof:}}
\def\endd{$\hfill\Box$}  				

\newcommand{\rmm}[1]{
	{\rm #1}
}

\newcommand{\tela}[1]{
	\tilde{#1}
}

\newcommand{\rmnum}[1]{\romannumeral #1}
\newcommand{\Rmnum}[1]{\uppercase\expandafter{\romannumeral#1}}

\newlist{rlist}{enumerate}{1}
\setlist[rlist]{label=\rmm{(}\!\!\roman*\!\rmm{)}}

\begin{document}
\begin{center}
{\Large  \bf  Intersection of complete cotorsion pairs}

\vspace{0.5cm} Qikai Wang and Haiyan Zhu\footnote{Corresponding author

\ Supported by the National Natural Science Foundation of China  (12271481).}  \\
School of Mathematical Science, Zhejiang University of Technology, Hangzhou 310023, China

E-mails: qkwang@zjut.edu.cn and hyzhu@zjut.edu.cn
\end{center}

\bigskip
\centerline { \bf  Abstract}
\leftskip10true mm \rightskip10true mm \noindent


Given two complete cotorsion pairs $(\mathcal{X}_1,\mathcal{Y}_1)$ and $(\mathcal{X}_2,\mathcal{Y}_2)$ in an exact category with $\mathcal{X}_1\subseteq \mathcal{Y}_2$,
we prove that $\left({\rm Smd}\langle \mathcal{X}_1,\mathcal{X}_2 \rangle,\mathcal{Y}_1\cap \mathcal{Y}_2\right)$ is also a complete cotorsion pair, where ${\rm Smd}\langle \mathcal{X}_1,\mathcal{X}_2 \rangle$ is the class of direct summands of extension of $\mathcal{X}_1$ and $\mathcal{X}_2$.
As an application, we  construct complete  cotorsion pairs,
such as $(^\perp\mathcal{GI}^{\leqslant n},\mathcal{GI}^{\leqslant n})$,  where $\mathcal{GI}^{\leqslant n}$ is  the class of modules of Gorenstein injective dimension at most $n$.
And we also characterize the  left orthogonal class of exact complexes of injective modules and the classes of modules with finite Gorenstein projective, Gorenstein flat, and PGF dimensions.
\\[2mm]

\leftskip0true mm \rightskip0true mm

\setcounter{tocdepth}{1} 

\include*{Introduction}

\include*{Main}

\include*{Exam}

\bibliography{tex}

\vspace{4mm}

\end{document}

%% file: Introduction.tex
\section {Introduction}

Ever since Salce introduced the notion of a cotorsion pair in the late 1970’s in \cite{Salce}, the significance of complete cotorsion pairs has been widely understood in approximation theory \cite{GR06}.
It is known that the intersection of the left (respectively, right)-hand classes of two cotorsion pairs is the left (respectively, right)-hand class of a cotorsion pair (See Enochs \cite[Exercise 7.1.2]{EE}).
In general, the intersection of the left(right)-hand class of complete cotorsion pairs is not the left(right)-hand class of a complete cotorsion pair, such as, \cite[Example]{HD05}.
Meanwhile, some examples of complete cotorsion pair in \cite[Table 1]{Zhang} with intersection are given by Gao, Lu and Zhang.
Consequently, one may inquire: Under what conditions can we get a new complete cotorsion pair  from the intersection of  two cotorsion pairs in exact categories.

The concept of an exact category is due to D. Quillen \cite{Q-1973}.
An \emph{exact category} is an additive category $\mathcal{A}$ equipped with a class $\mathcal{E}$ of kernel-cokernel sequences $A \stackrel{s}{\dan}B\stackrel{t}{\man}C$ in $\mathcal{A}$ such that $s$ is the kernel of $t$ and $t$ is the cokernel of $s$.
The class $\mathcal{E}$ satisfies exact axioms, for details, we refer the reader to \cite[Definition 2.1]{BT2010}.
Given an exact category $\mathcal{A}$, we will call a kernel-cokernel sequence a \emph{conflation} if it is in $\mathcal{E}$.
The morphism $s$ in a conflation $A \stackrel{s}{\dan}B\stackrel{t}{\man}C$ is called an
\emph{inflation}
and the morphism $t$ is called a \emph{deflation}.

	A \emph{cotorsion pair} is a pair of classes of objects $\left(\mathcal{X},\mathcal{Y}\right)$ of an exact category $\mathcal{A}$ such that $\mathcal{X}^\perp=\mathcal{Y}$ and $\mathcal{X} =\ ^\perp \mathcal{Y}$.
	Here $\mathcal{X}^\perp$, the right orthogonal class of $\mathcal{X}$, is the class of objects $A\in\mathcal{A}$ such that $\rmm{Ext}^1(X,A)=0$ for all $X\in \mathcal{X}$, and similarly $^\perp \mathcal{Y}$, the left orthogonal class of $\mathcal{Y}$, is the class of objects $A\in \mathcal{A}$ such that $\rmm{Ext}^1(A,Y)=0$ for all $Y\in \mathcal{Y}$. Moreover,
	a cotorsion pair $\left(\mathcal{X},\mathcal{Y}\right)$ is \emph{complete} if for every $A\in \mathcal{A}$ there are conflations $ Y\dan X\man A $ and $ A\dan Y'\man X'$ with $X,X'\in \mathcal{X}$ and with $Y,Y'\in \mathcal{Y}$.
	A cotorsion pair $\left(\mathcal{X},\mathcal{Y}\right)$ in an exact category $\mathcal{A}$ is  \emph{hereditary} if $\mathcal{X}$ is closed under the kernel of deflations and $\mathcal{Y}$ is closed under the cokernel of inflations.
    In addition, in a weakly idempotent complete exact category $\mathcal{A}$,  a complete cotorsion pair $\left(\mathcal{X},\mathcal{Y}\right)$ is  \emph{hereditary} if and only if $\mathcal{X}$ is closed under the kernel of deflations or $\mathcal{Y}$ is closed under the cokernel of inflations. (See \v{S}\v{t}ov\'{\i}\v{c}ek \cite[Lemma 6.17]{j14}).

The main purpose of this note is to provide an answer  to the above question.
For two classes $\mathcal{X}_1, \mathcal{X}_2$ of objects in  an exact category $\mathcal{A}$, denote by $\left\langle\mathcal{X}_1,\mathcal{X}_2\right\rangle = \{X\mid \exists $ a conflation $  X_2\dan X\man X_1$ with $X_1\in \mathcal{X}_1, X_2\in \mathcal{X}_2 \}$ and by $\rmm{Smd}\left\langle\mathcal{X}_1,\mathcal{X}_2\right\rangle$ the class of direct summands of objects in $\left\langle\mathcal{X}_1,\mathcal{X}_2\right\rangle$.

\begin{theo}\label{main2}
   Let $\left(\mathcal{X}_1,\mathcal{Y}_1\right)$ and $\left(\mathcal{X}_2,\mathcal{Y}_2\right)$ be  complete cotorsion pairs in an exact category $\mathcal{A}$.
   Then we have the following conclusions:
   \begin{rlist}
       \item If $\mathcal{Y}_2\subseteq \mathcal{X}_1$, then $\left(\mathcal{X}_1\cap \mathcal{X}_2, \rmm{Smd}\left\langle\mathcal{Y}_1,\mathcal{Y}_2\right\rangle\right)$ is a complete cotorsion pair in $\mathcal{A}$.
       Furthermore, if  $(\mathcal{X}_1,\mathcal{Y}_1)$ and $(\mathcal{X}_2,\mathcal{Y}_2)$ are hereditary complete cotorsion pairs  and  $\mathcal{A}$ is  weakly idempotent complete, then $\left(\mathcal{X}_1\cap \mathcal{X}_2, \rmm{Smd}\left\langle\mathcal{Y}_1,\mathcal{Y}_2\right\rangle\right)$ is hereditary.
       \item If $\mathcal{X}_1\subseteq \mathcal{Y}_2$, then $\left(\rmm{Smd}\left\langle\mathcal{X}_1,\mathcal{X}_2\right\rangle, \mathcal{Y}_1\cap\mathcal{Y}_2\right)$ is a complete cotorsion pair in $\mathcal{A}$.
       Furthermore, if  $(\mathcal{X}_1,\mathcal{Y}_1)$ and $(\mathcal{X}_2,\mathcal{Y}_2)$ are hereditary complete cotorsion pairs  and  $\mathcal{A}$ is  weakly idempotent complete, then $\left(\rmm{Smd}\left\langle\mathcal{X}_1,\mathcal{X}_2\right\rangle, \mathcal{Y}_1\cap\mathcal{Y}_2\right)$ is hereditary.
   \end{rlist}
\end{theo}

In fact,  {\textacutedbl}\!\rmm{Smd}\!\!{\textacutedbl} cannot be dropped by Example \ref{Smd} and the converse of Theorem \ref{main2} is not true in general by Example \ref{condition}.

As applications, we get a lot of known and unkown cotorsion pairs (see Example \ref{exam}).
In particular, we show that $(^\perp\mathcal{GI}^{\leqslant n},\mathcal{GI}^{\leqslant n})$ is a hereditary complete cotorsion pair (see Corollary \ref{GIn}), where $\mathcal{GI}^{\leqslant n}$ is  the class of modules of Gorenstein injective dimension less than or equal to $n$.
Furthermore,  we obtain new characterizations of the left orthogonal class of exact complexes of injective modules  and the classes of modules with finite Gorenstein projective(Gorenstein flat, PGF) dimensions (see Corollary \ref{coro1}).

%% file: Main.tex
\section{Proof of Theorem \ref{main2}}\label{mainsec}

\subsection{}It is easy to prove that $\rmm{Ext}^1(X,Y)=0$, for any $X\in \mathcal{X}_1\cap \mathcal{X}_2$ and $Y\in \rmm{Smd}\langle\mathcal{Y}_1,\mathcal{Y}_2\rangle$.
It follows that  $(\mathcal{X}_1\cap \mathcal{X}_2)\subseteq {^\perp\rmm{Smd}\left\langle\mathcal{Y}_1,\mathcal{Y}_2\right\rangle}$ and $\rmm{Smd}\langle\mathcal{Y}_1,\mathcal{Y}_2\rangle\subseteq (\mathcal{X}_1\cap \mathcal{X}_2)^\perp$.

\subsection{}  $\mathcal{X}_1\cap\mathcal{X}_2$ is special precovering.

For any object $A\in\mathcal{A}$, there is a conflation $Y_1\dan X_1\dan A$ with $X_1\in\mathcal{X}_1$ and $Y_1\in \mathcal{Y}_1$.
In addition, there is also a conflation
$Y_2\dan X_2\man X_1$ with $X_2\in  \mathcal{X}_2$ and $Y_2\in \mathcal{Y}_2$.
Note that $Y_2\in \mathcal{Y}_2\subseteq \mathcal{X}_1$ and $\mathcal{X}_1$ is closed under extension.
Then $X_2\in \mathcal{X}_1$, that is, $X_2\in \mathcal{X}_1\cap \mathcal{X}_2$.
and therefore we have the following commutative diagram:
\begin{center}
	\begin{tikzcd}
		Y_2 \arrow[d,tail] \arrow[r, equal] & Y_2 \arrow[d,tail]   &   \\
		Y \arrow[d,two heads] \arrow[r,tail]   & X_2 \arrow[d,two heads] \arrow[r,two heads] & A \arrow[d, equal]  \\
		Y_1 \arrow[r,tail]    & X_1 \arrow[r,two heads]  & A
	\end{tikzcd}
\end{center}
where $Y$ is the pullback of $Y_1\ra X_1$ and $X_2\ra X_1$. It follows from 2.1 that $Y\in (\mathcal{X}_1\cap \mathcal{X}_2)^\perp$.
So we deduce that $\mathcal{X}_1\cap\mathcal{X}_2$ is special precovering from the conflation $Y\dan X_2\man A$.

\subsection{}   $\rmm{Smd}\langle\mathcal{Y}_1,\mathcal{Y}_2\rangle$ is special preenveloping.

For any object $A\in\mathcal{A}$, there is a conflation $A\dan Y_1\man X_1$ with $X_1\in \mathcal{X}_1$ and $Y_1\in \mathcal{Y}_1$.
In addition, there is also a conflation $Y_2\dan X_2\man X_1$ with $X_2\in\mathcal{X}_2$ and $Y_2\in \mathcal{Y}_2$.
Note that $Y_2\in \mathcal{Y}_2\subseteq \mathcal{X}_1$ and $\mathcal{X}_1$ is closed under extension.
Then $X_2\in \mathcal{X}_1$, that is, $X_2\in \mathcal{X}_1\cap \mathcal{X}_2$.
Thus we have the following commutative diagram:
\begin{center}
	\begin{tikzcd}
		&Y_2 \arrow[r,equal]\arrow[d,tail]&Y_2\arrow[d,tail]& \\
		A\arrow[d,equal] \arrow[r,tail]  & Y\arrow[r,two heads] \arrow[d,two heads]& X_2 \arrow[d,two heads] \\
		A\arrow[r,tail]& Y_1\arrow[r,two heads]&X_1
	\end{tikzcd}
\end{center}
where $Y$ is the pullback of $Y_1\ra X_1$ and $X_2\ra X_1$.
From the conflation $A\dan Y\man X_2$, we can deduce that $\rmm{Smd}\langle\mathcal{Y}_1,\mathcal{Y}_2\rangle$ is special preenveloping.

\subsection{} $^\perp\rmm{Smd}\langle\mathcal{Y}_1,\mathcal{Y}_2\rangle\subseteq (\mathcal{X}_1\cap \mathcal{X}_2)$.

For any object $A\in{^\perp\rmm{Smd}\langle\mathcal{Y}_1,\mathcal{Y}_2\rangle}$, there is a conflation  $Y\rightarrowtail X\twoheadrightarrow A$    with $X \in \mathcal{X}_1\cap \mathcal{X}_2$ and $Y\in\rmm{Smd}\langle\mathcal{Y}_1,\mathcal{Y}_2\rangle$ by 2.2.
Note that $\rmm{Ext}^1(A,Y)=0$ by assumption, and hence the conflation $Y\dan X\man A$ is split.
Thus $A$ is a  direct summand  of $X$, which implies   $A\in \mathcal{X}_1\cap \mathcal{X}_2$.

\subsection{ }  $(\mathcal{X}_1\cap\mathcal{X}_2)^\perp\subseteq \rmm{Smd}\langle\mathcal{Y}_1,\mathcal{Y}_2\rangle$.

For any object $A\in (\mathcal{X}_1\cap\mathcal{X}_2)^\perp$, there is a conflation $A\dan Y\man X$ with $Y\in\langle\mathcal{Y}_1,\mathcal{Y}_2\rangle$ and $X\in\mathcal{X}_1\cap\mathcal{X}_2$ by 2.3.
Since  $\rmm{Ext}^1(X,A)=0$ by assumption, that is, $Y\simeq A\bigoplus X$ which  induces  $A\in \rmm{Smd}\langle\mathcal{Y}_1,\mathcal{Y}_2\rangle$.

Therefore,
$(\mathcal{X}_1\cap \mathcal{X}_2, \rmm{Smd}\langle\mathcal{Y}_1,\mathcal{Y}_2\rangle)$ is a complete cotorsion pair in $\mathcal{A}$.
And, it is easy to prove that $\mathcal{X}_1\cap \mathcal{X}_2$ is closed under the kernel of deflations if  $(\mathcal{X}_1,\mathcal{Y}_1)$ and $(\mathcal{X}_2,\mathcal{Y}_2)$ are hereditary complete cotorsion pairs.
\endd

\begin{remark}
    If the category $\mathcal{A}$ have enough projective objects, then we only need to show that 2.1, 2.3 and 2.4 by Salce’s Lemma. Dually,  If the category $\mathcal{A}$ have enough projective objects or enough injective objects, then we only need to show that 2.1, 2.2 and 2.5.
\end{remark}

%% file: Exam.tex
\section{Examples}

To better apply the theorem, we need the following definitions for  a ring $R$ with identity.
From now on, and unless otherwise stated, by a module we mean a left $R$-module.

Recall that  a complex is \emph{exact} if it is exact at all terms.
And a complex $P_\bullet$ is \emph{dg-projective} if $P_n$ is a projective module for every $n\in \mathbb{Z}$ and  $\rmm{Hom}(P_\bullet,E_\bullet)$ is exact for any exact complex $E_\bullet$.
For details, we refer the reader to \cite[Chapter 4]{EE2}.

A module  is called \emph{Gorenstein projective}  \cite{AM1969} if it is a cycle  of an exact complex of projective modules which remains exact after applying $\rmm{Hom}_R(-,P)$ for any projective module $P$.

A module  is  called \emph{Gorenstein flat}  \cite{EE} (\emph{projectively coresolved Gorenstein flat} \cite{SS}) if it is a cycle  of an exact complex  of flat (projective) modules which remains exact after applying $E_R\otimes-$ for any  right injective module $E$.

The $\rmm{FP}_n$-modules were introduced (\cite{BDPM}) as a refinement of the class of finitely presented modules.
The $\rmm{FP}_n$-injective modules are the modules in the right orthogonal class of $\rmm{FP}_n$-modules.
We define the \emph{Gorenstein $\rmm{FP}_n$-injective} modules to be the cycles of the exact complexes of injective modules which remains exact after apply  $\rmm{Hom}(A, -)$ for any $\rmm{FP}_n$-injective module $A$.
Given that $\rmm{FP}_0$-injective modules are equivalent to injective modules, it follows that Gorenstein $\rmm{FP}_0$-injective modules are simply Gorenstein injective modules.

Let $\mathcal{X}$ be a class  of  modules.
We say that the $\mathcal{X}$-dim$M\leqslant n$ if there exists an exact complex
$$0\ra X_{n}\ra \cdots\ra X_1\ra X_0\ra M\ra 0$$
with $X_0,X_1,\dots,X_{n-1} \in \mathcal{X}$.
Dually, we have the definition of  $\mathcal{X}$-codim$M$.

For convenience, we give the list of notations.

\begin{table}[htbp]
    \centering
    \begin{tabular}{ll}
        \toprule[1pt]
        \hspace{0.3cm} Symbol&\hspace{1.5cm} Description(for any positive integer $n$)\\
        \midrule[1pt]
        \hspace{0.5cm} $\mathcal{P}$ & \hspace{1.5cm} the class of projective $R$-modules.  \vspace{-3pt}\\
        \hspace{0.5cm} $\mathcal{I}$ & \hspace{1.5cm} the class of injective $R$-modules.  \vspace{-3pt}\\
        \hspace{0.5cm} $\mathcal{F}$ & \hspace{1.5cm} the class of flat $R$-modules. \vspace{-3pt}\\
        \hspace{0.5cm} $\mathcal{P}^{\leqslant n}$ &\hspace{1.5cm} the class of $R$-modules $M$ with   $\mathcal{P}$-dim$M\leqslant n$. \vspace{-3pt}\\
        \hspace{0.5cm} $\mathcal{I}^{\leqslant n}$ &\hspace{1.5cm} the class of $R$-modules $M$ with   $\mathcal{I}$-codim$M\leqslant n$. \vspace{-3pt}\\
        \hspace{0.5cm} $\mathcal{F}^{\leqslant n}$&\hspace{1.5cm} the class of $R$-modules $M$ with  $\mathcal{F}$-dim$M\leqslant n$. \vspace{-3pt}\\
        \hspace{0.5cm} $\mathcal{GP}$ &\hspace{1.5cm} the class of Gorenstein projective  modules.\vspace{-3pt}\\
        \hspace{0.5cm} $\mathcal{GP}^{\leqslant n}$ &\hspace{1.5cm} the class of $R$-modules $M$ with  $\mathcal{GP}$-dim$M\leqslant n$. \vspace{-3pt}\\
        \hspace{0.5cm} $\mathcal{GF}$&\hspace{1.5cm} the class of  Gorenstein flat modules.\vspace{-3pt}\\
        \hspace{0.5cm} $\mathcal{GF}^{\leqslant n}$&\hspace{1.5cm} the class of $R$-modules $M$ with  $\mathcal{GF}$-dim$M\leqslant n$. \vspace{-3pt}\\
        \hspace{0.5cm} $\mathcal{PGF}$&\hspace{1.5cm} the class of  $\rmm{PGF}$ modules.\vspace{-3pt}\\
        \hspace{0.5cm} $\mathcal{PGF}^{\leqslant n}$&\hspace{1.5cm} the class of $R$-modules  $M$ with  $\mathcal{PGF}$-dim$M\leqslant n$. \vspace{-3pt}\\
     
        \hspace{0.5cm} $\mathcal{FI}_n$&\hspace{1.5cm} the class of  $\rmm{FP}_n$-injective modules.\vspace{-3pt}\\

        \hspace{0.5cm} $\mathcal{FI}_1^{\leqslant n}$&\hspace{1.5cm} the class of modules of  $M$ with  $\mathcal{FI}_1$-codim$M\leqslant n$.		\vspace{-3pt}\\
        \hspace{0.5cm} $\mathcal{GI}_n$&\hspace{1.5cm} the class of Gorenstein $\rmm{FP}_n$-injective modules.\vspace{-3pt}\\
        \hspace{0.5cm} $\mathcal{GI}$&\hspace{1.5cm} the class of Gorenstein injective modules.\vspace{-3pt}\\
        \hspace{0.5cm} $\rmm{dg}\tela{\mathcal{P}}$  &\hspace{1.5cm} the class of dg-projective complexes.\vspace{-3pt}\\
        \hspace{0.5cm} $\mathcal{E}$&\hspace{1.5cm} the class of  exact  complexes.\vspace{-3pt}\\
        \hspace{0.5cm} $\rmm{dw}\tela{\mathcal{I}}$   &\hspace{1.5cm} the class of  complexes $X_\bullet$ with every $X_n$  injective.\vspace{-3pt}\\
        \hspace{0.5cm} $\rmm{ex}\tela{\mathcal{I}}$ &\hspace{1.5cm} the class of  exact  complexes $X_\bullet$ with every $X_n$  injective.\\
        \bottomrule[1pt]
    \end{tabular}
\end{table}

To obtain a new cotorsion pair, we need the following inclusion relation.

\begin{lemm}\label{thick}\label{new2}
    For any positive integer $n$ and any integer $m\geqslant 2$, $\mathcal{FI}_1^{\leqslant n}\subseteq{^\perp \mathcal{GI}_m}$.
\end{lemm}
\proof
By the definitions of $\mathcal{GI}_m$, $\mathcal{FI}_n$ and \cite[Lemma 2.2]{GIn} ,we have the following diagram:
\begin{center}
    \begin{tikzcd}
        \mathcal{FI}_0 \arrow[r, "\subseteq"] \arrow[d, "\subseteq"] & \mathcal{FI}_1 \arrow[r, "\subseteq"] \arrow[d, "\subseteq"] & \dots \arrow[d, "\subseteq"] \arrow[r, "\subseteq"] & \mathcal{FI}_n \arrow[d, "\subseteq"] \\
        ^\perp\mathcal{GI} \arrow[r, "\subseteq"]                    & ^\perp\mathcal{GI}_1 \arrow[r, "\subseteq"]                  & \dots \arrow[r, "\subseteq"]                        & ^\perp\mathcal{GI}_n
    \end{tikzcd}

    Diag. 1
\end{center}
Let $M\in\mathcal{FI}_1^{\leqslant n}$.
Then there is an exact complex
\begin{center}
    $0\ra M\ra E_0\ra E_1\ra \cdots \ra E_{n-1} \ra E_n\ra 0$.
\end{center}
with $E_i\in \mathcal{FI}_1$.
Let  $K_i$  be the kernel of the morphism  $E_{i} \ra E_{i+1}$, for $0<i\leqslant n$.
Consider the short exact sequence
$0\ra K_{n-1}\ra E_{n-1}\ra E_{n}\ra 0$.
Since $E_{n-1},E_{n}\in \mathcal{FI}_1\subseteq {^\perp\mathcal{GI}_m}$,  the thick class $^\perp \mathcal{GI}_m$ $($\cite[Corollary 2.10]{GIn}$)$ give us that $K_{n-1}\in{^\perp\mathcal{GI}_m}$.
Similarly, the short exact sequences $0\ra K_{n-2}\ra E_{n-2}\ra K_{n-1}\ra0$ show that $k_{n-2}\in{^\perp\mathcal{GI}_m}$.
By induction, $M\in{^\perp\mathcal{GI}_m}$, that is, $\mathcal{FI}_1^{\leqslant n}\subseteq{^\perp \mathcal{GI}_m}$.
\endd

\begin{exam}\label{exam}
   The complete cotorsion pairs in the right column of the following table are obtained by Theorem \ref{main2}, where  $n> 0$ and $m\geqslant2$.

   \begin{table}[htbp]
       \setlength{\tabcolsep}{0pt}
       \centering

       \begin{tabular}{lll}
           \toprule[1pt]
            {\small Complete Cotorsion Pairs\hspace{0.6cm} $+$}\hspace{0.6cm} & {\small Inclusions\hspace{1.3cm} $\Longrightarrow$}\hspace{0.5cm} &{\small Induced Complete Cotorsion Pairs}\hspace{-30pt}\\
           \midrule[1pt]
           \hspace{0.5cm}{\small $\left( \mathcal{P}^{\leqslant n}, (\mathcal{P}^{\leqslant n})^\perp\right)$}&\vspace{-7pt}&\\
           \hspace{0.5cm}{\tiny\cite[Theorem 7.4.6]{EE}}& &\vspace{-1pt}\\
           \hspace{0.5cm}{\small $\left( \mathcal{GP},\mathcal{GP}^\perp\right)$}& {\small $\mathcal{P}^{\leqslant n}\!\subseteq\! \mathcal{GP}^\perp$ }  &{\small $\left(\rmm{Smd}\langle \mathcal{P}^{\leqslant n},\mathcal{GP}\rangle, (\mathcal{P}^{\leqslant n})^\perp\!\cap\! \mathcal{GP}^\perp\right)$}\vspace{-7pt} \\
           \hspace{0.5cm}\tiny{$($Artin algebra$)$} & {\tiny\cite[Theorem 2.5]{Zhang}} &\tiny{$($Artin algebra$)$} \vspace{-7pt}\\
           \hspace{0.5cm}{\tiny\cite[\Rmnum{10}, Theorem 2.4(\rmnum{4})]{BAR} }&  & \vspace{-1pt}\\
           \hspace{0.5cm}{\small $\left(\mathcal{F}^{\leqslant n},(\mathcal{F}^{\leqslant n})^\perp\right)$}&{\small  $\mathcal{P}^{\leqslant n}\!\subseteq\!\mathcal{F}^{\leqslant n}\!\subseteq\! \mathcal{PGF}^\perp$}&{\small $\left( \rmm{Smd}\langle \mathcal{P}^{\leqslant n},\mathcal{PGF}\rangle, (\mathcal{P}^{\leqslant n})^\perp\!\cap\! \mathcal{PGF}^\perp\right)$}\vspace{-7pt}\\
           \hspace{0.5cm}{\tiny\cite[Theorem 3.4(2)]{MD07}}&{\tiny\cite[Theorem 3.4]{REM2023} } &\vspace{-1pt}\\

           \hspace{0.5cm}{\small $\left(\mathcal{PGF},\mathcal{PGF}^\perp\right)$} &  & {\small $\left(\rmm{Smd}\langle\mathcal{F}^{\leqslant n},\mathcal{PGF}\rangle,(\mathcal{F}^{\leqslant n})^\perp\cap \mathcal{PGF}^\perp\right)$} \vspace{-7pt} \\
           \hspace{0.5cm}{\tiny\cite[Theorem 4.9]{SS}}&&\\
           \midrule[0.2pt]
           \hspace{0.5cm}{\small  $\left( ^\perp\mathcal{I}^{\leqslant n},\mathcal{I}^{\leqslant n}\right)$}&\small{$\mathcal{I}^{\leqslant n}\!\subseteq\! {^\perp \mathcal{GI}}$} &\vspace{-7pt} \\
           \hspace{0.5cm}{\tiny\cite[Theorem 8.7]{RJ06}}&{\tiny\cite[Corollary 11.2.7]{EE}} &\vspace{-1pt}\\
           \hspace{0.5cm}{\small  $\left( ^\perp\mathcal{GI},\mathcal{GI}\right)$}&\small{$\mathcal{FI}_m\!\subseteq\! {^\perp \mathcal{GI}_m}$} & $\left( {^\perp \mathcal{GI}}\!\cap\!{^\perp \mathcal{I}^{\leqslant n}}, \rmm{Smd}\langle  \mathcal{GI}, \mathcal{I}^{\leqslant n}\rangle \right)$\vspace{-7pt} \\
           \hspace{0.5cm}{\tiny\cite[Theorem 5.6]{SS}}&{\tiny see Diag. 1}  &\vspace{-1pt}\\
           \hspace{0.5cm}{\small $\left({^\perp\mathcal{FI}_n},\mathcal{FI}_n\right)$}&  {\small $\mathcal{FI}_n\!\subseteq\!\mathcal{FI}_{n+1}$}& {\small $\left(^\perp \mathcal{GI}\!\cap\!{^\perp \mathcal{FI}_0}, \rmm{Smd}\langle\mathcal{GI},\mathcal{FI}_0\rangle\right)$}\vspace{-7pt}\\
           \hspace{0.5cm}{\tiny\cite[Corollary 4.2]{BDPM} }&{\tiny by definition} &\vspace{-1pt}\\
           \hspace{0.5cm}{\small $\left({^\perp\mathcal{GI}_m},\mathcal{GI}_m\right)$} & {\small  ${^\perp \mathcal{GI}_m}\!\subseteq\! {^\perp \mathcal{GI}_{m+1}}$}& {\small $\left(^\perp \mathcal{GI}_m\!\cap\!{^\perp \mathcal{FI}_n}, \rmm{Smd}\langle\mathcal{GI}_m,\mathcal{FI}_n\rangle\right)$}\vspace{-7pt} \\
           \hspace{0.5cm}{\tiny\cite[Theorem 2.24]{GIn}}&{\tiny by definition}& {\tiny  $(n\leqslant m)$\vspace{-1pt}}\\
           \hspace{0.5cm}{\small $\left( ^\perp\mathcal{FI}^{\leqslant n}_1,\mathcal{FI}^{\leqslant n}_1\right)$} &$\mathcal{FI}_1^{\leqslant n}\!\subseteq\!{^\perp \mathcal{GI}_m}$ & {\small $\left({^\perp \mathcal{GI}_m}\!\cap\!{^\perp\mathcal{FI}^{\leqslant n}_1},\rmm{Smd}\langle\mathcal{GI}_m, \mathcal{FI}^{\leqslant n}_1\rangle\right)$}\vspace{-7pt}\\
           \hspace{0.5cm}\tiny{$($coherent ring$)$} &{\tiny By Lemma \ref{new2} } &{\tiny $($coherent ring$)$}\vspace{-7pt}\\
           \hspace{0.5cm}{\tiny\cite[Theorem 3.8]{MD05}} &  &\\
           \midrule[0.2pt]

           \hspace{0.5cm}{\small$\left(\rmm{dg}\tela{\mathcal{P}},\mathcal{E}\right)$} & \vspace{-7pt}\\
           \hspace{0.5cm}{\tiny\cite[Chapter 4]{EE2}}& &\vspace{-1pt}\\
           \hspace{0.5cm}{\small $\left(^\perp \rmm{dw}\tela{\mathcal{I}},\rmm{dw}\tela{\mathcal{I}}\right)$} & {\small $^\perp \rmm{dw}\tela{\mathcal{I}}\!\subseteq\!  \mathcal{E}$}&{\small $\left(\rmm{Smd}\langle^\perp \rmm{dw}\tela{\mathcal{I}},\rmm{dg}\tela{\mathcal{P}}\rangle, \rmm{dw}\tela{\mathcal{I}}\cap\mathcal{E}\right)$}\vspace{-7pt}\\
           \hspace{0.5cm}{\tiny\cite[Proposition 3.2]{dw}}&{\tiny by definition} &\\
           \bottomrule[1pt]
       \end{tabular}
   \end{table}
\end{exam}

Next  we give the new characterizations of $R$-modules with finite Gorenstein projective dimension, Gorenstein flat dimension, $\rmm{PGF}$ dimension and the left orthogonal class of $\rmm{ex}\tela{\mathcal{I}}$.
\vspace{-5pt}
\begin{coro}\label{coro1}
    For any ring $R$ and every $n>0$,
    we have the following conclusions:
    \begin{rlist}
        \item $ \mathcal{PGF}^{\leqslant n}= \rmm{Smd}\langle \mathcal{P}^{\leqslant n},\mathcal{PGF}\rangle$
        \item $ \mathcal{GF}^{\leqslant n}= \rmm{Smd}\langle\mathcal{F}^{\leqslant n},\mathcal{PGF}\rangle$
        \item  $^\perp \rmm{ex}\tela{\mathcal{I}} = \rmm{Smd}\langle{^\perp \rmm{dw}\tela{\mathcal{I}}},\rmm{dg}\tela{\mathcal{P}}\rangle$
        \item if $R$ is an Artin algebra, then $\mathcal{GP}^{\leqslant n}= \rmm{Smd}\langle \mathcal{P}^{\leqslant n},\mathcal{GP}\rangle$.
    \end{rlist}
\end{coro}
\vspace{-8pt}
\proof
By \cite{Zhang} and Example \ref{exam}, (\emph{\rmnum{1}}),(\emph{\rmnum{2}}) and (\emph{\rmnum{4}}) hold.
In addition, (\emph{\rmnum{3}}) holds by \cite{dw} and Example \ref{exam}.
\endd
\vspace{-5pt}

\begin{coro}\label{GIn}
    For any ring and $n\geqslant0$, $({^\perp \mathcal{GI}^{\leqslant n}}, \mathcal{GI}^{\leqslant n})$ is a hereditary complete cotorsion pair.
\end{coro}
\vspace{-8pt}
\proof
By the Example \ref{exam}, we only have to show that $\mathcal{GI}^{\leqslant n}= \rmm{Smd}\langle \mathcal{GI},\mathcal{I}^{\leqslant n}\rangle$.

First, we show that $\mathcal{GI}^{\leqslant n}\subseteq \rmm{Smd}\langle \mathcal{GI},\mathcal{I}^{\leqslant n}\rangle= ({^\perp \mathcal{GI}}\!\cap\!{^\perp \mathcal{I}^{\leqslant n}})^\perp$.
For any module $M\in \mathcal{GI}^{\leqslant n}$, there is a short exact sequence $0\ra M\ra G\ra N\ra 0$ with $G\in \mathcal{GI}$ and $N\in \mathcal{I}^{\leqslant n-1}$ by \cite[ Theorem 2.15]{HH2004}.
Since $G\in \mathcal{GI}$, there is a short exact sequence $0\ra G'\ra I\ra G\ra 0$ with $I\in \mathcal{I}$ and $G'\in \mathcal{GI}$ by definition.
Then, we have the following commutative diagram:
\begin{center}
    \begin{tikzcd}
         G' \arrow[d,tail] \arrow[r,equal] & G' \arrow[d,tail]          &                          \\
         L \arrow[d,two heads] \arrow[r,tail]  & I \arrow[d,two heads] \arrow[r,two heads] & N \arrow[d,equal]  \\
         M \arrow[r,tail]   & G \arrow[r,two heads]  & N
    \end{tikzcd}
\end{center}
where $L$ is the pullback of $M\ra G$ and $I\ra G$.
Then $L\in \mathcal{I}^{\leqslant n}$.
For any $X\in {^\perp\mathcal{GI}}\cap{^\perp \mathcal{I}^{\leqslant n}} $,
we have the following commutative diagram with exact rows
\begin{center}
    \begin{tikzcd}
     {\rmm{Hom}(X,N)} \arrow[d,equal] \arrow[r] & {\rmm{Ext}^1(X,L)}\arrow[d] &\\
     {\rmm{Hom}(X,N)} \arrow[r]                    &   {\rmm{Ext}^1(X,M)} \arrow[r] & {\rmm{Ext}^1(X,G)}
    \end{tikzcd}
\end{center}
By ${\rmm{Ext}^1(X,L)}=0={\rmm{Ext}^1(X,G)}$,
it follows that ${\rmm{Ext}^1(X,M)}=0$.
So $M\in ({^\perp\mathcal{GI}}\cap{^\perp \mathcal{I}^{\leqslant n}})^\perp$, as desired.

Conversely, for any module $M\in \rmm{Smd}\langle \mathcal{GI}, \mathcal{I}^{\leqslant n} \rangle$,
there is an exact sequence
$ 0\ra N\ra M\bigoplus L \ra G\ra 0$
with $G\in \mathcal{GI}$ and $N\in \mathcal{I}^{\leqslant n}$.
Clearly, the Gorenstein injective dimensions of $G$ and $N$ are finite, and so is $M$ by \cite[Theorem 2.24]{HH2004}.
For any $I\in \mathcal{I}$ and $i>n$, we have the following exact sequence
$$\rmm{Ext}^i(I,N)\ra\rmm{Ext}^i(I,M\bigoplus L)\ra\rmm{Ext}^i(I,G)$$
Since $I\in \mathcal{I}\subseteq {^\perp\mathcal{GI}}$,  $\rmm{Ext}^i(I,G)=0$.
And $N\in \mathcal{I}^{\leqslant n}$ implies $\rmm{Ext}^i(I,N)=0$.
It follows that $\rmm{Ext}^i(I,M\bigoplus L)=0$, that is, $\rmm{Ext}^i(I,M)=0$ for all $i>n$.
Therefore, $M\in \mathcal{GI}^{\leqslant n}$ by \cite[Theorem 2.22]{HH2004}.
\endd

\vspace{15pt}
The following example from \cite[Page 15]{LW00} shows that the {\textgravedbl\!\rmm{Smd}\!\!\textacutedbl}  cannot be dropped.
\begin{exam}\label{Smd}
    Let $k$ be a field, and $R=k\llbracket X,Y\rrbracket/(XY)$.
    Then we claim $^\perp((\mathcal{P}^{\leqslant 1})^\perp\!\cap\! \mathcal{GP}^\perp) =\mathcal{GP}^{\leqslant 1}\neq \langle\mathcal{P}^{\leqslant 1},\mathcal{GP} \rangle$.
    Note that $R$ is a $1$-dimensional Gorenstein ring but not a domain and, hence, not regular.
    Then $k$ is a module in $\mathcal{GP}^{\leqslant 1}$.
    Assume that there is a short exact sequence
    \begin{center}
          $ 0\ra G \ra k\ra P_1\ra 0$
    \end{center}
    with $ G\in \mathcal{GP}$ and $P_1\in \mathcal{P}_1$.
    Since $k$ is a simple module,
    $k\cong GP$ or
    $k\cong P_1$.
    However, the fact that $R$  is not regular induces that $Pd_Rk=\infty$.
    And $k$ is not Gorenstein projective,
    a contradiction.
    Therefore, $k$ is not a extension of $P_1$ and $GP$.
\end{exam}

Next,   the converse of Theorem \ref{main2}(\emph{\rmnum{2}}) is not true by the following example.

\begin{exam}\label{condition}
    Consider the finite dimensional $K$-algebra $\mathcal{A}=KQ/I$, where $Q$ is the quiver
    \begin{center}
        \begin{tikzcd}
            4 \arrow[r, "\alpha"] \arrow[d, "\alpha'"] & 3 \arrow[d, "\beta"] \\
            2 \arrow[r, "\beta'"]                      & 1
        \end{tikzcd}
    \end{center}
    and $I=\langle\beta\alpha - \beta'\alpha'\rangle$.
    The category $\mathcal{A}$-mod of finite dimensional left $\mathcal{A}$-modules, has eleven indecomposable modules,
    \begin{center}
        \rmm{ind}$\mathcal{A}$-mod $= \{P_1, P_2, P_3, P_4, I_2,I_3, I_4\simeq S_4, S_2, S_3, M_1\simeq (P_2\oplus P_3)/ P_1, M_2 \simeq (S_2\oplus S_3\oplus P_4)/M_1 \}$
    \end{center}
    It is well-known that $\mathcal{A}$-mod has the Auslander-Reiten quiver
    \begin{center}
        \begin{tikzcd}
            & P_2 \arrow[rd] \arrow[rr, no head, dotted] &                                     & S_3 \arrow[rd] \arrow[rr, no head, dotted] &                                                       & I_2 \arrow[rd] &     \\
            P_1 \arrow[ru] \arrow[rd] \arrow[rr, no head, dotted] &                                            & M_1 \arrow[ru] \arrow[rd] \arrow[r] & P_4 \arrow[r]                              & M_2 \arrow[ru] \arrow[rd] \arrow[rr, no head, dotted] &                & I_4 \\
            & P_3 \arrow[ru] \arrow[rr, no head, dotted] &                                     & S_2 \arrow[ru] \arrow[rr, no head, dotted] &                                                       & I_3 \arrow[ru] &
        \end{tikzcd}
    \end{center}
    According to the Auslander-Reiten quiver, we have the following almost split sequences:
    \begin{align*}
        P_1 \ra P_2\ra S_2 & & P_3\ra M_1\ra S_2& & P_3\ra P_4 \ra I_2  \\
        M_1\ra P_4\oplus S_2\ra I_2 & & S_3\ra M_2\ra I_2& &
    \end{align*}
    Set
    \begin{align*}
        \mathcal{D}=\{P_1,P_2,P_3,P_4,I_2,M_1,S_2\},&&
        \mathcal{D}_1=\{P_1,P_2,P_3,P_4,I_2\} &\ \ and&
        \mathcal{D}_2=\{P_1,P_2,P_3,P_4,M_1,S_2\}.
    \end{align*}
    Then
    \begin{align*}
        \mathcal{D}^\perp&=\{ P_2, P_4, I_2, I_3, I_4, M_2, S_2 \} \\
        \mathcal{D}_1^\perp&=\{ P_1, P_2, P_4, I_2, I_3, I_4, M_2, S_2 \}\\
        \mathcal{D}_2^\perp&=\{ P_2, P_4, I_2, I_3, I_4, M_1, M_2, S_2, S_3\}.
    \end{align*}
    It is easy to check that $(\mathcal{D},\mathcal{D}^\perp), (\mathcal{D}_1,\mathcal{D}_1^\perp), (\mathcal{D}_2,\mathcal{D}_2^\perp)$ are cotorsion pairs. And they are complete by above almost split sequences and Salce's Lemma.
    But, $\mathcal{D}_1\subsetneq \mathcal{D}_2^\perp$ and $\mathcal{D}_2\subsetneq \mathcal{D}_1^\perp$ since $P_3\notin \mathcal{D}_1^\perp\cup \mathcal{D}_2^\perp$.
\end{exam}

\section*{Acknowledgements}

The authors would like to thank Jiangsheng Hu, Xiaoyan Yang  and Sergio Estrada for interesting discussions.
The authors are also grateful to the referees for the valuable comments and suggestions.

%% file: tex.bib
@incollection {j14,
    AUTHOR = {\v{S}\v{t}ov\'{\i}\v{c}ek, J.},
     TITLE = {Exact model categories, approximation theory, and cohomology
              of quasi-coherent sheaves},
 BOOKTITLE = {Advances in representation theory of algebras},
    SERIES = {EMS Ser. Congr. Rep.},
     PAGES = {297--367},
 PUBLISHER = {Eur. Math. Soc., Z\"{u}rich},
      YEAR = {2014},
   MRCLASS = {18E10 (18E30 18F20)},
  MRNUMBER = {3220541},
MRREVIEWER = {R. H. Street},
}

@book {RJ06,
    AUTHOR = {G\"{o}bel, R. and Trlifaj, J.},
    TITLE = {Approximations and endomorphism algebras of modules. {V}olume
    2},
    PUBLISHER = {Walter de Gruyter GmbH \& Co. KG, Berlin},
    YEAR = {2012},
    PAGES = {i--xxiv and 459--972},
    ISBN = {978-3-11-021810-7; 978-3-11-021811-4},
    MRCLASS = {16-02 (03E75 16Dxx 16S50)},
    MRNUMBER = {2985654},
    MRREVIEWER = {J. L. G\'{o}mez Pardo},
    DOI = {10.1515/9783110218114},
    URL = {https://doi.org/10.1515/9783110218114},
}

@book {LW00,
    AUTHOR = {Christensen, L. W.},
    TITLE = {Gorenstein dimensions},
    PUBLISHER = {Springer-Verlag, Berlin},
    YEAR = {2000},
    PAGES = {viii+204},
    ISBN = {3-540-41132-1},
    MRCLASS = {13D05 (13D25)},
    MRNUMBER = {1799866},
    MRREVIEWER = {Peter J\o rgensen},
    DOI = {10.1007/BFb0103980},
    URL = {https://doi.org/10.1007/BFb0103980},
}

@book {GR06,
    AUTHOR = {G\"{o}bel, R. and Trlifaj, J.},
    TITLE = {Approximations and endomorphism algebras of modules},
    PUBLISHER = {Walter de Gruyter GmbH \& Co. KG, Berlin},
    YEAR = {2006},
    PAGES = {xxiv+640},
    ISBN = {978-3-11-011079-1; 3-11-011079-2},
    MRCLASS = {16D90 (03E75 16D10 16E30 16S50 16S90)},
    MRNUMBER = {2251271},
    MRREVIEWER = {Sverre O. Smal\o },
    DOI = {10.1515/9783110199727},
    URL = {https://doi.org/10.1515/9783110199727},
}

@article {HD05,
    AUTHOR = {Happel, D. and Unger, L.},
    TITLE = {On a partial order of tilting modules},
    JOURNAL = {Algebr. Represent. Theory},
    FJOURNAL = {Algebras and Representation Theory},
    VOLUME = {8},
    YEAR = {2005},
    NUMBER = {2},
    PAGES = {147--156},
    ISSN = {1386-923X},
    MRCLASS = {16G10 (16E10 16G70)},
    MRNUMBER = {2162278},
    MRREVIEWER = {Apostolos D. Beligiannis},
    DOI = {10.1007/s10468-005-3595-2},
    URL = {https://doi.org/10.1007/s10468-005-3595-2},
}

@article {MD07,
    AUTHOR = {Mao, L.X. and Ding, N.Q.},
    TITLE = {Envelopes and covers by modules of finite {FP}-injective and
    flat dimensions},
    JOURNAL = {Comm. Algebra},
    FJOURNAL = {Communications in Algebra},
    VOLUME = {35},
    YEAR = {2007},
    NUMBER = {3},
    PAGES = {833--849},
    ISSN = {0092-7872},
    MRCLASS = {16E10 (16D40 16D50)},
    MRNUMBER = {2305235},
    MRREVIEWER = {Javad Asadollahi},
    DOI = {10.1080/00927870601115757},
    URL = {https://doi.org/10.1080/00927870601115757},
}

@article {MD05,
    AUTHOR = {Mao, L.X. and Ding, N.Q.},
    TITLE = {Relative {FP}-projective modules},
    JOURNAL = {Comm. Algebra},
    FJOURNAL = {Communications in Algebra},
    VOLUME = {33},
    YEAR = {2005},
    NUMBER = {5},
    PAGES = {1587--1602},
    ISSN = {0092-7872},
    MRCLASS = {16E10 (16D40)},
    MRNUMBER = {2149078},
    MRREVIEWER = {Thomas Cassidy},
    DOI = {10.1081/AGB-200061047},
    URL = {https://doi.org/10.1081/AGB-200061047},
}

@incollection {Salce,
    AUTHOR = {Salce, L.},
    TITLE = {Cotorsion theories for abelian groups},
    BOOKTITLE = {Symposia {M}athematica, {V}ol. {XXIII} ({C}onf. {A}belian
    {G}roups and their {R}elationship to the {T}heory of
    {M}odules, {INDAM}, {R}ome, 1977)},
    PAGES = {11--32},
    PUBLISHER = {Academic Press, London-New York},
    YEAR = {1979},
    MRCLASS = {20K40 (18E40)},
    MRNUMBER = {565595},
    MRREVIEWER = {P. L. Sperry},
}

@article {dw,
    AUTHOR = {Gillespie, J.},
    TITLE = {Cotorsion pairs and degreewise homological model structures},
    JOURNAL = {Homology Homotopy Appl.},
    FJOURNAL = {Homology, Homotopy and Applications},
    VOLUME = {10},
    YEAR = {2008},
    NUMBER = {1},
    PAGES = {283--304},
    ISSN = {1532-0073},
    MRCLASS = {18G55 (55U35)},
    MRNUMBER = {2399475},
    MRREVIEWER = {David A. Blanc},
    DOI = {10.4310/HHA.2008.v10.n1.a12},
    URL = {https://doi.org/10.4310/HHA.2008.v10.n1.a12},
}

@article {BDPM,
    AUTHOR = {Bravo, D. and P\'{e}rez, M. A.},
    TITLE = {Finiteness conditions and cotorsion pairs},
    JOURNAL = {J. Pure Appl. Algebra},
    FJOURNAL = {Journal of Pure and Applied Algebra},
    VOLUME = {221},
    YEAR = {2017},
    NUMBER = {6},
    PAGES = {1249--1267},
    ISSN = {0022-4049},
    MRCLASS = {16E05 (16E30 16S90 18G25)},
    MRNUMBER = {3599429},
    MRREVIEWER = {David A. Jorgensen},
    DOI = {10.1016/j.jpaa.2016.09.008},
    URL = {https://doi.org/10.1016/j.jpaa.2016.09.008},
}

@article {BT2010,
    AUTHOR = {B\"{u}hler, T.},
    TITLE = {Exact categories},
    JOURNAL = {Expo. Math.},
    FJOURNAL = {Expositiones Mathematicae},
    VOLUME = {28},
    YEAR = {2010},
    NUMBER = {1},
    PAGES = {1--69},
    ISSN = {0723-0869},
    MRCLASS = {18E10 (18-02 18E30)},
    MRNUMBER = {2606234},
    MRREVIEWER = {Sunil K. Chebolu},
    DOI = {10.1016/j.exmath.2009.04.004},
    URL = {https://doi.org/10.1016/j.exmath.2009.04.004},
}

@inproceedings {Q-1973,
    AUTHOR = {Quillen, D.},
    TITLE = {Higher algebraic {$K$}-theory. {I}},
    BOOKTITLE = {Algebraic {$K$}-theory, {I}: {H}igher {$K$}-theories ({P}roc.
    {C}onf., {B}attelle {M}emorial {I}nst., {S}eattle, {W}ash.,
    1972)},
    SERIES = {Lecture Notes in Math., Vol. 341},
    PAGES = {85--147},
    PUBLISHER = {Springer, Berlin-New York},
    YEAR = {1973},
    MRCLASS = {18F25},
    MRNUMBER = {338129},
    MRREVIEWER = {Stephen M. Gersten},
}

@article {GIn,
    AUTHOR = {Iacob, A.},
    TITLE = {Generalized {G}orenstein modules},
    JOURNAL = {Algebra Colloq.},
    FJOURNAL = {Algebra Colloquium},
    VOLUME = {29},
    YEAR = {2022},
    NUMBER = {4},
    PAGES = {651--662},
    ISSN = {1005-3867},
    MRCLASS = {18G10 (16D40 16D50 16D90 18G25 18G35)},
    MRNUMBER = {4517666},
    MRREVIEWER = {R. Hafezi},
    DOI = {10.1142/S1005386722000463},
    URL = {https://doi.org/10.1142/S1005386722000463},
}

@article {REM2023,
    AUTHOR = {El Maaouy, R.},
     TITLE = {Model structures, {$n$}-{G}orenstein flat modules and {PGF}
              dimensions},
   JOURNAL = {Proc. Edinb. Math. Soc. (2)},
  FJOURNAL = {Proceedings of the Edinburgh Mathematical Society. Series II},
    VOLUME = {67},
      YEAR = {2024},
    NUMBER = {4},
     PAGES = {1241--1264},
      ISSN = {0013-0915},
   MRCLASS = {42B25 (16E10 18G65 18G80 18N40)},
  MRNUMBER = {4832996},
       DOI = {10.1017/S0013091524000518},
       URL = {https://doi.org/10.1017/S0013091524000518},
}

@misc{Zhang,
    title={Chains of model structures arising from modules of finite Gorenstein dimension, arXiv 2403.05232 (math.{RT})},
    author={Gao, N.  and  Lu, X.S. and  Zhang, P.},
    eprint={2403.05232},
    archivePrefix={arXiv},
    primaryClass={math.RT}
}

@book {EE,
    AUTHOR = {Enochs, E. E. and Jenda, O. M. G.},
    TITLE = {Relative homological algebra},
    PUBLISHER = {Walter de Gruyter \& Co., Berlin},
    YEAR = {2000},
    ISBN = {3-11-016633-X},
    MRCLASS = {16E65 (13H10 16-02 16E10 18G25)},
    MRNUMBER = {1753146},
    MRREVIEWER = {J. Kuzmanovich},
    DOI = {10.1515/9783110803662},
    URL = {https://doi.org/10.1515/9783110803662},
}

@book {EE2,
    AUTHOR = {Enochs, E. E. and Jenda, O. M. G.},
    TITLE = {Relative homological algebra. {V}olume 2},
    PUBLISHER = {Walter de Gruyter GmbH \& Co. KG, Berlin},
    YEAR = {2011},
    ISBN = {978-3-11-021522-9},
    MRCLASS = {18G25 (16-02 16E05 18-02)},
    MRNUMBER = {2841617},
    MRREVIEWER = {Sergio Estrada-Dominguez},
}

@article {BAR,
    AUTHOR = {Beligiannis, A. and Reiten, I.},
    TITLE = {Homological and homotopical aspects of torsion theories},
    JOURNAL = {Mem. Amer. Math. Soc.},
    FJOURNAL = {Memoirs of the American Mathematical Society},
    VOLUME = {188},
    YEAR = {2007},
    NUMBER = {883},
    ISSN = {0065-9266},
    MRCLASS = {18G55 (20J05 55U35)},
    MRNUMBER = {2327478},
    DOI = {10.1090/memo/0883},
    URL = {https://doi.org/10.1090/memo/0883},
}

@article {SS,
    AUTHOR = {\v{S}aroch, J. and \v{S}\v{t}ov\'{\i}\v{c}ek, J.},
    TITLE = {Singular compactness and definability for {$\Sigma$}-cotorsion
    and {G}orenstein modules},
    JOURNAL = {Selecta Math. (N.S.)},
    FJOURNAL = {Selecta Mathematica. New Series},
    VOLUME = {26},
    YEAR = {2020},
    NUMBER = {2},
    PAGES = {Paper No. 23, 40},
    ISSN = {1022-1824},
    MRCLASS = {16E30 (03E75 16B70)},
    MRNUMBER = {4076700},
    MRREVIEWER = {Ramalingam Udhayakumar},
    DOI = {10.1007/s00029-020-0543-2},
    URL = {https://doi.org/10.1007/s00029-020-0543-2},
}

@book {AM1969,
    AUTHOR = {Auslander, M. and Bridger, M.},
    TITLE = {Stable module theory},
    PUBLISHER = {American Mathematical Society, Providence, RI},
    YEAR = {1969},
    PAGES = {146},
    MRCLASS = {16.40 (18.00)},
    MRNUMBER = {269685},
    MRREVIEWER = {G. Michler},
}

@article {HH2004,
    AUTHOR = {Holm, H.},
    TITLE = {Gorenstein homological dimensions},
    JOURNAL = {J. Pure Appl. Algebra},
    FJOURNAL = {Journal of Pure and Applied Algebra},
    VOLUME = {189},
    YEAR = {2004},
    NUMBER = {1--3},
    PAGES = {167--193},
    ISSN = {0022-4049},
    MRCLASS = {16E10 (16E05 16E30)},
    MRNUMBER = {2038564},
    MRREVIEWER = {Zhaoyong Huang},
    DOI = {10.1016/j.jpaa.2003.11.007},
    URL = {https://doi.org/10.1016/j.jpaa.2003.11.007},
}
